\documentclass[12pt]{amsart}

\usepackage{graphicx,amscd,hyperref,tikz-cd,mathtools,float}
\usepackage{enumitem}

\numberwithin{equation}{section}

 \newtheorem{theorem}{Theorem}[section]
    \newtheorem{prop}[theorem]{Proposition}

    \newtheorem{rem}[theorem]{Remark}
    \newtheorem{definition}[theorem]{Definition}
    \newtheorem{cor}[theorem]{Corollary}

\newcommand{\ov}{\overline}

\newcommand{\fm}{\mathfrak{m}}

\newcommand{\fp}{\mathfrak{p}}

\newcommand{\G}{\textnormal{SL}_2\mathbb{C}}

\newcommand{\Ke}{K_\epsilon}
\newcommand{\Kz}{K_\zeta}

\newcommand{\be}{\begin{equation}}
  \newcommand{\ee}{\end{equation}}

\newcommand{\Spec}{\mathrm{Spec}}
\newcommand{\MS}{\mathrm{Max} \ \mathrm{Spec}}
\newcommand{\End}{\mathrm{End}}

\usepackage{pict2e}
\makeatletter
\newcommand*\@KP@Large@frame[2]{%
    \setlength\unitlength{\fontdimen 22 #1\tw@}%
    \vrule \@width\z@ \@height 4\unitlength \@depth\tw@\unitlength
    \begin{picture}(6,2)(-3,-1)%
        \def\@KP@Radius     {3}%
        \def\@KP@Hole@radius{.5}% The same value seem adequate for both...
        \def\@KP@Diameter   {6}%
        #2%
    \end{picture}%
}
\newcommand*\@KP@Small@frame[2]{%
    \setlength\unitlength{\fontdimen 22 #1\tw@}%
    \vrule \@width\z@ \@height \thr@@\unitlength \@depth\@ne\unitlength
    \begin{picture}(4,2)(-2,-1)%
        \def\@KP@Radius     {2}%
        \def\@KP@Hole@radius{.5}% ... but let it be customizable too.
        \def\@KP@Diameter   {4}%
        #2%
    \end{picture}%
}

\newcommand*\@KP@Radius     {}
\newcommand*\@KP@Hole@radius{}
\newcommand*\@KP@Diameter   {}
\newcommand*\@KP@Shape@A{%
    \put(0,0){\circle{\@KP@Diameter}}%
}
\newcommand*\@KP@Shape@B{%
    \Line(-\@KP@Radius,\@KP@Radius )(\@KP@Radius,-\@KP@Radius)%
    \Line(-\@KP@Radius,-\@KP@Radius)(-\@KP@Hole@radius,-\@KP@Hole@radius)%
    \Line(\@KP@Radius ,\@KP@Radius )(\@KP@Hole@radius ,\@KP@Hole@radius )%
}
\newcommand*\@KP@Shape@E{%
    \Line(-\@KP@Radius,\@KP@Radius )(\@KP@Radius,-\@KP@Radius)%
     \Line(\@KP@Radius,\@KP@Radius )(-\@KP@Radius,-\@KP@Radius)%

}

\newcommand*\@KP@Shape@C{%
    \cbezier(-\@KP@Radius,\@KP@Radius )(0,0)(0,0)(\@KP@Radius,\@KP@Radius )%
    \cbezier(-\@KP@Radius,-\@KP@Radius)(0,0)(0,0)(\@KP@Radius,-\@KP@Radius)%
}
\newcommand*\@KP@Shape@D{%
    \cbezier(-\@KP@Radius,-\@KP@Radius)(0,0)(0,0)(-\@KP@Radius,\@KP@Radius)%
    \cbezier(\@KP@Radius ,-\@KP@Radius)(0,0)(0,0)(\@KP@Radius ,\@KP@Radius)%
}

\newcommand*\@KP@Atomic@mathpalette[1]{%
    \mathinner{% or "\mathord"?
        \mathchoice{%
            \linethickness{.6\p@}% Tip: use thicker lines if you decide to 
            \@KP@Large@frame \textfont {#1}%
        }{%
            \linethickness{.4\p@}% adjustable
            \@KP@Small@frame \textfont {#1}%
        }{%
            \linethickness{.3\p@}% adjustable
            \@KP@Small@frame \scriptfont {#1}%
        }{%
            \linethickness{.2\p@}% adjustable
            \@KP@Small@frame \scriptscriptfont {#1}%
        }%
    }%
}

\newcommand*\KPA{\@KP@Atomic@mathpalette \@KP@Shape@A}
\newcommand*\KPB{\@KP@Atomic@mathpalette \@KP@Shape@B}
\newcommand*\KPC{\@KP@Atomic@mathpalette \@KP@Shape@C}
\newcommand*\KPD{\@KP@Atomic@mathpalette \@KP@Shape@D}
\newcommand*\KPE{\@KP@Atomic@mathpalette \@KP@Shape@E}

\title[Skein Module at an irreducible representation] {The Kauffman Bracket Skein Module\\ at an irreducible representation}
\author[F. Tehrani]{Mohammad Farajzadeh-Tehrani}\author[Frohman]{Charles Frohman} \author[Kania-Bartoszynska]{Joanna Kania-Bartoszynska}

\begin{document}
\maketitle
\begin{abstract} In this paper, we study the Kauffman bracket skein module of closed oriented three-manifolds at a non-multiple-of-four roots of unity.  Our main result establishes that the localization of this module at a maximal ideal, which corresponds to an irreducible representation of the fundamental group of the manifold, forms a one-dimensional free module over the localized unreduced coordinate ring of the character variety. 
We apply this by proving that the dimension of the skein module of a homology sphere with finite character variety, when the order of the root of unity is not divisible by $4$,  is greater than or equal to the dimension of the unreduced coordinate ring of the character variety. This leads to a computation of the dimension of the skein module with coefficients in rational functions for homology spheres with  tame universal skein module.
 \end{abstract}

\tableofcontents
\section{Introduction}

The Kauffman bracket skein module of an oriented three-manifold $X$ at a root of unity $\zeta$, denoted   $\Kz(X)$, is formed  by taking the complex vector space with basis the isotopy classes of framed links in $X$ and modding out by the submodule corresponding to Kauffman bracket skein relations with parameter $\zeta$. 
If $\zeta$ is a primitive $n$-th root of unity where $n$ is not divisible by $4$, let 
$$
m=\frac{n}{\gcd(n,2)}\qquad \mathrm{and}\qquad \epsilon=\zeta^{m^2}\in \{\pm 1\}.
$$ 
The Kauffman bracket skein module $\Ke(X)$ is defined similarly, only using $\epsilon$ as the variable in the Kauffman bracket skein relations. The module $\Ke(X)$  is an algebra over the complex numbers under disjoint union and it isomorphic to the unreduced coordinate ring of the $\mathrm{SL}(2,{\mathbb C} )$-character variety of the fundamental group of the manifold $X$; see \cite{Bu,PS}.\\

The threading map of Bonahon and Wong \cite{BW,Le1} makes $\Kz(X)$ into a module over $\Ke(X)$. This paper explores $\Kz(X)$ from the viewpoint of this module structure with the goal of understanding quantum invariants of framed links in three-manifolds. \\

The Kauffman bracket skein module was inspired by the existence of  the Reshetikhin-Turaev invariant \cite{Li}. It is constructed so that the Reshetikhin-Turaev invariant 
$$ 
RT^\zeta:\Kz(X)\rightarrow \mathbb{C}
$$ 
is a linear functional. For a long time, this was the only linear functional on the Kauffman bracket skein module that was not the result of an ad-hoc construction. Later, the modified trace \cite{CGP} introduced more linear functionals.  The recent paper \cite{FKL}  showed that for each irreducible representation of the fundamental group there is a linear functional on $\Kz(X)$.  These linear functionals come from specializing $\Kz(X)$ at the maximal ideal of $\Ke(X)$ associated to irreducible representations of $\pi_1(X)$.  This led to advances in the study of the Kauffman bracket skein module \cite{DKS,KK}. \\
  
  The Kauffman bracket skein module of $X$ defined over $\mathbb{Q}(q)$ is being studied with an eye towards the construction of topological field theories \cite{GJS}. It is denoted by $K_q(X)$. Detcherry, Kalfagianni, and Sikora \cite{DKS} prove that the dimension of $\Kz(X)$ as a complex vector space is greater than or equal to the number of maximal ideals of $\Ke(X)$ in the case that $X$ has a finite character variety. They used this to estimate the dimension of $K_q(X)$,  with the added assumption that the universal skein module of $X$ is tame, in terms of the dimension of $\Kz(X)$.   Leveraging standard constructions in the study of modules over an affine commutative ring, we are able to strengthen these results in the case that $X$ is a homology sphere. \\

  {\bf Theorem \ref{itsbig}. }  {\em Let $X$ be an oriented homology sphere with finite character variety. If $\zeta$ is a primitive root of unity of order $n$ where $n$ is not divisible by $4$,  then the dimension of $\Kz(X)$ as a complex vector space is greater than or equal to the dimension of $\Ke(X)$.}

 \begin{rem}
 A closed three-manifold is {\bf small} if it has no incompressible surfaces. It is a Theorem of Culler and Shalen that if a three-manifold is small then its $\G$-character variety is finite \cite{CS}. For instance, most three-manifolds obtained by surgery on a hyperbolic knot are small.  
\end{rem}

The technical advance of this paper is the following result, which also applies to three-manifolds whose character variety has positive dimensional-components.\\

{\bf Theorem \ref{main_thm}.} {\em
Suppose that $\zeta$ is a primitive root of unity of order $n$ where $n$ is not divisible by $4$,  and let $\epsilon=\zeta^{(\frac{n}{ \gcd(n,2)})^2}.$ 
Let $X$ be a closed oriented $3$-manifold, and  $\fm\leq \Ke(X)$  the maximal ideal corresponding to an irreducible representation $\rho:\pi_1(X)\rightarrow \G$. The localized skein module $\Kz(X)_\fm $  is a free module of rank one over $\Ke(X)_\fm$; i.e. $\Kz(X)_\fm\cong \Ke(X)_\fm$. }\\

Furthermore, in the context of Theorem~\ref{itsbig}, if the universal Kauffman bracket skein module $K(X)$ is tame and has no torsion coming from the cyclotomic polynomial associated to $\zeta$ then $\Kz(X)$ is a free module of rank one over $\Ke(X)$ and we get the following result. \\
%Here $K_q(X)$ denotes the Kauffman bracket skein module with coefficients in rational functions $\mathbb{Q}(q)$. %Already defined above.

{\bf Theorem \ref{rational}.} {\em  If $X$ is a homology sphere with  tame universal Kauffman bracket skein module, then 
$$  
\dim_{\mathbb{C}}\Ke(X)=\dim_{\mathbb{Q}(q)}K_q(X).
$$ }

From the work of \cite{DKS} we know that if the universal skein module is tame then the character variety of $\pi_1(X)$ is finite. It should be noted that if the character variety of $\pi_1(X)$ is infinite, nothing like this can happen since $K_q(X)$ is  finite dimensional \cite{GJS}.  \\

The modern vision of skein modules is that they come from tangle functors defined in balls via some sort of technique for extension to three-manifolds \cite{GJS, CGP1}. An upshot of the explorations in this paper is that the construction of tangle functors is an unfinished business.  There should be tangle functors where there is a background flat connection on the ball  that takes into account infinitesimal deformations of the flat connection. \\

The paper is organized as follows. In Section~\ref{revkbsm}, we review skein modules and skein algebras. In Section~\ref{Algebra_sec}, we review the definitions around localization. We also review the properties of left and right ideals of a matrix algebra and their quotients. In Section~\ref{Skein_sec}, we discuss the computation of the skein module of a closed three-manifold at a root of unity and prove preparatory lemmas needed for the proof of the main result. After proving some statements about local modules at an irreducible representation in Section~\ref{localized-model_sec}, the proof of Theorem~\ref{main_thm} will be completed in  Section~\ref{proof_sec}. Finally, in Section~\ref{cons} we prove that the dimension $\Kz(X)$ when $X$ is a homology sphere is greater than or equal to the dimension of $\Ke(X)$.   \\

This paper was inspired by the recent paper \cite{DKS} that used properties of skein modules at roots of unity
to draw conclusions about skein modules defined over the rational functions. The authors are deeply appreciative of Frauke Bleher for reading and editing early versions of this paper.  We also  thank Renaud Detcherry, Effie Kalfagianni and Adam Sikora for their input on an earlier draft of this paper.

\section{Review of the Kauffman bracket skein modules and algebras}\label{revkbsm} 

In this section we recall the necessary definitions and known results about the structure of the Kauffman bracket skein module of an oriented $3$-manifold and its algebra structure when the $3$-manifold is a cylinder over a surface.
 
\begin{definition} \label{kbsm} Let $X$ be a closed oriented $3$-manifold and $\mathcal{L}(X)$ denote the set of isotopy classes of framed links in $X$, including the empty link. The {\bf universal Kauffman bracket skein module} $K(X)$ is the quotient of the free  $\mathbb{Z}[q,q^{-1}]$-module (here, $q$ is a formal variable) with basis $\mathcal{L}(X)$, by its submodule of  Kauffman bracket skein relations
$$
\aligned
&\KPA + q^{2}+q^{-2}=0 \\
  & \KPB -q \KPC - q^{-1}\KPD=0.
        \endaligned
$$ 
If $R$ is any ring and $\zeta\in R$ is a unit, then $R$ is a module over $\mathbb{Z}[q,q^{-1}]$ where $q$ acts as multiplication by $\zeta$. The {\bf Kauffman bracket skein module with variable $\zeta$} is the specialization of $K(X)$,
$$
K_\zeta(X)=K(X)\otimes_{\mathbb{Z}[q,q^{-1}]} R.
$$

\end{definition}

In the last section we will study
$$ 
K_q(X)=K(X)\otimes_{\mathbb{Z}[q,q^{-1}]}\mathbb{Q}(q) 
$$ 
where $\mathbb{Q}(q)$ is the field of rational functions in the variable $q$. 
However, we are primarily interested in the cases where $R=\mathbb{C}$ and $\zeta$ is a  root of unity whose order  $n$ is not divisible by $4$. In this case, let 
\be\label{mep_e}
m=\frac{n}{\textnormal{gcd}(n,2)}\quad\mathrm{and}\quad \epsilon=\zeta^{m^2}\in \{\pm 1\}.
\ee
The skein module $\Ke(X)$ is an algebra. Multiplication on $\Ke(X)$ comes from taking the disjoint union of embedded framed links. 
\begin{theorem}\label{Dougsbigthm}\cite{Bu}, \cite{PS1}.
The algebra $\Ke(X)$ is isomorphic to the unreduced coordinate ring of the $\G$-character variety of $\pi_1(X)$  
\end{theorem}

In passing to $\Kz(X)$, one loses the algebra structure but gains the ability to measure the entanglement of framed links. \\

For $\zeta$ and $\epsilon$ as in (\ref{mep_e}), the  {\bf threading} map of Bonahan and Wong  \cite{BW,Le1}
\be\label{tau}
\tau\equiv \tau_m:K_\epsilon(X)\rightarrow K_\zeta(X),
\ee
comes from threading links with a Chebyshev polynomial. The Chebyshev polynomials  $\{T_k\}_{k\in \mathbb{N}}$ of the first type  are defined recursively by  
\be \label{ChP-FT_e}
T_0(x) =2,\quad T_1(x)= x,~\mathrm{and}\quad T_{k+1}(x) = x T_k(x) - T_{k-1}(x)\quad \forall~k\geq 1.
\ee
The threading map descends from the linear map
$$
T\colon \mathbb{C}\mathcal{L}(X)\rightarrow \mathbb{C}\mathcal{L}(X),
$$
which systematically replaces each component of a framed link in a multilinear fashion by applying $T_m(x)$ to that component, guided by the framing-annulus associated with that particular component.\\

The threading map is natural in the sense that it commutes with inclusions.  Also, the image of any framed link under the threading map is {\bf transparent}; the skein represented by its  disjoint union with any other link is independent of crossing changes with the threaded link. Therefore, the threading map realizes $\Kz(X)$ as a  module over $\Ke(X)$.  The $\Ke(X)$-module structure on $\Kz(X)$ comes from threading followed by taking disjoint union.  We will see later that $\Kz(X)$ is a finitely generated module over $\Ke(X)$.\\

If $X=F\times [0,1]$ is a cylinder over a surface, we denote the skein module by $\Kz(F)$. In this case, $\Kz(F)$ is an algebra under the multiplication given by stacking framed links on  top of each other. It is well-known that  $\Kz(F)$ is a finitely generated algebra (c.f. \cite{B}) without zero divisors (c.f. \cite{PS}).  The threading map  
$
\tau\colon \Ke(F)\rightarrow \Kz(F) 
$ 
is injective and its image is contained in the center $Z(\Kz(F))$.  The skein algebra  $\Kz(F)$ is  a finite rank module over $\Ke(F)$; see \cite{AF}.\\

Let $\Sigma_{g,b}$ denote the compact oriented surface of genus $g$ with $b$ boundary components.  The peripheral skeins are the skeins coming from blackboard framed simple closed curves on the surface $\Sigma_{g,b}$ that are parallel to the boundary components. Let $\partial_i$ denote the peripheral skein that is parallel to the $i$-th boundary component. 
\begin{theorem}\label{center}\cite{FKL1},\cite {FKL2}
Let $\zeta$ be a root of unity of order not divisible by $4$ and $m$ be the integer defined in (\ref{mep_e}). The center of the Kauffman bracket skein algebra of $\Sigma_{g,b}$,  when $\chi(\Sigma_{g,b})<0$, 
 is the polynomial algebra
$$
Z(\Kz(\Sigma_{g,b}))=\tau(\Ke(\Sigma_{g,b}))[\partial_1,\ldots,\partial_b].
$$
Let $S= Z(\Kz(\Sigma_{g,b}))-\{0\}$. The dimension of the localization $S^{-1}\Kz(\Sigma_{g,b})$ as a vector space over $S^{-1}Z(\Kz(\Sigma_{g,b}))$ is equal $m^{6g-6+2b}$.  
\end{theorem}

\section{Review of the algebraic foundations}\label{Algebra_sec}
There are two related operations on a module $M$ over a commutative algebra $A$. If $\mathfrak{a}\leq A$ is an ideal, the {\bf specialization} of $M$ at $\mathfrak{a}$ is
\be 
M/\mathfrak{a}M=M\otimes_A A/\mathfrak{a}.
\ee  
If $S\subset A$ is a multiplicatively closed subset of $A$, then the {\bf localization} $S^{-1}A$ of $A$ with respect to $S$ is an algebra obtained by inverting the elements of $S$ and the localization of  $M$ with respect to $S$ is the $S^{-1}A$-module  
\be
M\otimes_A S^{-1}A.
\ee  
We are particularly interested in the case where $\fm \leq A$ is a maximal ideal and $S=A-\fm$. In this case, we denote the localizations of $A$ and $M$ with respect to $A-\fm$ by $A_\fm$ and $M_\fm$, respectively.\\

The goal of this section is to state the algebraic results needed to relate the operations of specialization and localization  of the Kauffman bracket skein module of a closed oriented $3$-manifold.

\subsection{Localization, factorization and Artinian rings} 
Let $A$ be a commutative algebra over the complex numbers. The {\bf spectrum} of $A$, denoted $\Spec(A)$, is the set of prime ideals of $A$.  If $s\in A$, then the principal open set determined by $s$ is 
\begin{equation}\label{D(s)_eq}
D(s)=\{\fp \in \Spec(A)|s\not\in \fp\}.
\ee
The principal open sets form a basis for the {\bf Zariski topology}  on $\Spec(A)$.  In the Zariski topology, the closure of a prime ideal $\fp$ is the set of all prime ideals containing $\fp$. The set of maximal ideals of $A$, denoted $\MS(A)$, thus corresponds to the set of closed points in $\Spec(A)$.\\

Every  algebra morphism $h:A\rightarrow B$ induces a continuous map
$$
h^*\colon\Spec(B)\rightarrow \Spec(A)
$$ 
that sends the prime ideal $\fp\in \Spec(B)$ to  $h^*(\fp)=h^{-1}(\fp)$.   If $h$ is surjective, then
$h^*$ embeds $\Spec(B)$ as a closed subset of $\Spec(A)$. Specifically,
$$ 
h^*(\Spec(B))=\{\fp \in \Spec(A)| \ker(h)\leq \fp\}.
$$ 
To minimize the proliferation of symbols in the context of an algebra $B$ that is a surjective image of $A$, we opt to identify $\Spec(B)$ with its image in $\Spec(A)$. This approach allows us to name ideals only once, streamlining the notation.\\

For a non-necessarily commutative algebra $A$, let $Z(A)$ denote its center.   
Given a multiplicatively closed subset $S\subset Z(A)$, the {\bf localization}
$S^{-1}A$ of $A$ with respect to $S$ is the set of equivalence classes of ordered pairs $(a,s)\in A\times S$ subject to the relation 
$$
(a,s)\sim (a',s') \Leftrightarrow \exists u\in S \quad \textnormal{s.t.}\quad u(as'-a's)=0.
$$  
The equivalence classes are multiplied and added as if they were fractions with denominators in $S$, to induce a canonical algebra structure on $S^{-1}A$. Thus, it is expedient to denote the equivalence class of $(a,s)$ by $a/s$. It is worth noting that the natural algebra homomorphism $A\rightarrow S^{-1}A$ given by $a\rightarrow a/1$ is not injective when $A$ has zero divisors.\\

If $M$ is a left or right module over $A$, then the localization of $M$ at $S$, denoted $S^{-1}M$, is a module over $S^{-1}A$. We restrict our comments to left modules, the formulas for right modules are similar. Let  $S^{-1}M$  be the set of equivalence classes of ordered pairs $(m,s)\in M\times S$ subject to the relation 
$$
(m,s)\sim (m',s') \Leftrightarrow \exists u\in S \quad \textnormal{s.t.}\quad u(s'm-sm')=0.
$$  
 Here juxtaposition is used to denote the action of $A$ on $M$. Equivalence classes of  ordered pairs are added as if they are fractions with denominator in $S$ and $S^{-1}A$ acts on $S^{-1}M$ as fractions. Alternatively, one can carry out the localization of a module  by  taking the tensor product. If $M$ is a left module over $A$, then
$$ 
S^{-1}M=S^{-1}A\otimes_AM.
$$ 
Localization is an exact functor that commutes with sums, quotients and tensor products. If $S$ is the set of powers of a single
non-nilpotent element $s\in A$, then we denote the localizations of $A$ and $M$ by $A_s$ and $M_s$, respectively.  In the case where $S$ is the complement of a prime ideal $\fp$, the localizations are denoted by $A_\fp$ and $M_\fp$, respectively.\\

Recall a partially ordered set is {\bf directed} if for every $i,j\in I$ there exists $k\in I$ with $i,j\leq k$. A directed system of modules $(M_i,\mu_{ij})$ is a collection of modules over an algebra $A$  indexed by a directed set $I$ and morphisms $\mu_{ij}:M_i\rightarrow M_j$ for each   $i\leq j$. The morphisms have  the properties that if $i<j<k$ then $\mu_{jk}\circ \mu_{ij}=\mu_{ik}$ and $\mu_{ii}=Id:M_i\rightarrow M_i$.  Let $C$ be the direct sum of the $M_i$ where we identify each $M_i$ with its image in $C$. Let $D\leq C$ be the submodule spanned by all $\mu_{ij}(x_i)-x_i$ where $x_i\in M_i$. The direct limit is defined to be
\be \varinjlim_{i\in I}M_i=C/D.\ee  Inclusion followed by the quotient map yields
\be \mu_i:M_i\rightarrow \varinjlim_{i\in I} M_i,\ee from which we can see that every
$x\in \varinjlim_{i\in I}M_i$ can be represented by
$\mu_i(x_i)$
for some $i\in I$. Also if $\mu_i(x_i)=0$ then $\mu_{ij}(\mu_i(x_i))=0$ for some $j>i$, see \cite{AM}.
\\

Suppose that $A$ is a commutative algebra. With notation as in (\ref{D(s)_eq}), if $\{D(s)\}_{s\in I}$ is a neighborhood basis of the maximal ideal $\fm\in \Spec(A)$, then the collection $\{A_s\}_{s\in I}$ can be seen as  a directed system of modules over $A$ that are also algebras. Since $D(t)\!\subset\! D(s)$ if and only if $t$ is in the radical of the principal ideal generated by $s$, there is $b\in A$ and a natural number $n$ so that $t^n=bs$. The morphisms 
\be \label{directed} 
A_s\rightarrow A_t,\qquad a/s^l\rightarrow ab^l /t^{ln}\qquad \forall~a\in A,l\in \mathbb{N},
\ee
are independent of $b$ and $n$, and make the $\{A_s\}_{s\in I}$ into a directed system of modules over $A$ satisfying
$$
A_\fm=\varinjlim_{s\in I} A_s. 
$$
Since the morphisms in the directed system are also algebra morphisms the limit,  $A_\fm$ is an algebra.
In the language of  sheaves, the stalk of the structure sheaf of $\Spec(A)$ at $\fm$ is $A_{\fm}$; c.f. \cite{GW}.
If $M$ is a module over $A$, since $M_s=M \otimes_A A_s$ and $M_\fp=M\otimes_A A_\fp$, we similarly get
\be\label{dirmod}
M_\fm=\varinjlim_{s\in I} M_s.
\ee

When $M$ is a module over a commutative algebra $A$, and
 if 
\be \label{A-decomposition_e}
A\cong \prod_{i=1}^nA/\mathfrak{a}_i
\ee
such that $\mathfrak{a_i}\!\leq\! A$ are ideals, then the module $M$ also factors as a product 
\be \label{M-decomposition_e}
M\cong \prod_{i=1}^nM/\mathfrak{a}_iM,
\ee
c.f. \cite{AM}. In the context of Theorems~\ref{itsbig} and~\ref{rational}, we are particularly interested in the case where $\mathfrak{a}_1= \fm^{k}$ and $\fm$ is  the maximal ideal corresponding to an isolated (possibly fat, that is with nonzero Zariski tangent space) point   and $\mathfrak{a}\equiv \mathfrak{a}_2$ is the ideal corresponding to the complement of that point. In this case
\be \label{Max-factor_e}
A\cong (A/\fm^k)\times A/\mathfrak{a}\quad \textnormal{and} \quad M\cong (M/\fm^kM)\times M/\mathfrak{a}
\ee
for some $k\in\mathbb{Z}_+$. 
In our specific example of interest in this paper, the positive integer $k$ will be  the order of contact between two holomorphic Lagrangians at that  point.\\

 An algebra is deemed {\bf Artinian} if every descending chain of ideals stabilizes; c.f. \cite{AM}. Every commutative algebra over $\mathbb{C}$ that has  finite dimension as a vector space is Artinian. If $A$ is Artinian, it has a finite number of maximal ideals, denoted  $\{\fm_i\}_{i=1}^n$, and there exist positive integers $k_i$ such that $A \cong \prod_{i=1}^nA/(\fm_i)^{k_i}.$  
Consequently, for any module $M$ over $A$, the isomorphism $M \cong \prod_{i=1}^n \left(M/\fm_i^{k_i}M\right)$ holds true. In the special case of a three-manifold with finite character variety, the Kauffman bracket skein algebra $\Ke(X)$ is Artinian.\\

In (\ref{A-decomposition_e}), if some factor $A/\mathfrak{a}_i$ is a local Artinian ring, then the corresponding factor $M/\mathfrak{a}_iM$ in (\ref{M-decomposition_e}) coincides with the localization of $M$ at the maximal ideal of the local Artinian ring $A/\mathfrak{a}_i$. In particular, the latter applies to the special situation of (\ref{Max-factor_e}): the specializations $A/\fm^{k}$ and $M/\fm^kM$ coincide with the localizations $A_{\fm}$ and $M_\fm$, respectively. 
\subsection{Left and right ideals of a  matrix algebra }  
In this paper, we  also deal with  matrix algebras having coefficients in a commutative algebra $A$ over $\mathbb{C}$.
More specifically, let $M_n(A)\cong A^n\otimes_AA^n$ be the algebra of $n\times n$ matrices with coefficients from $A$. 
If $L\leq M_n(A)$ is a left ideal, the {\bf row space}   $V(L)\leq A^n$ of $L$ is the span of all the rows of all the matrices in $L$. It is easy to see that 
 \begin{equation}\label{V(L)_e}
 A^n\otimes_AV(L)=L.
 \end{equation}   
 Similary if $R\leq M_n(A)$ is a right ideal, the column space  $V(R)\leq A^n$ of $R$ is the span of the columns of all the matrices in $R$.  We have  \begin{equation}\label{V(R)_e}
 V(R)\otimes_AA^n=R.
 \end{equation} 
 Left ideals are classified by their row space as a submodule of $A^n$, and right ideals are classified by their column space. \\

We use the following consequences of the right-exactness of tensor product, when working with left and right modules that arise from a Heegaard splitting.
If
$$
i_1\colon N_1\rightarrow  M_1\quad \textnormal{and}\quad  i_2\colon N_2\rightarrow M_2
$$ 
are morphisms of  $A$-modules then
\be \label{onside}
(M_1\otimes_A M_2)/(M_1\otimes_A i_2(N_2)) \cong M_1\otimes_A (M_2/i_2(N_2))
\ee
and
\be\label{MN-decomposition} 
\frac{M_1\otimes_A M_2}{M_1\otimes_A i_2(N_2)+i_1(N_1)\otimes_A M_2} \cong \frac{M_1}{i_1(N_1)}\otimes_A \frac{M_2}{i_2(N_2)}.
\ee

\begin{cor}\label{L+R_cor}
Suppose that  $L,R\leq M_n(A)$ are left and right ideals, respectively.  If  $V(L)\subset A^n$ is  the row space of $L$ and  $V(R)\subset A^n$ is  the column space of $R$,   then
  $$ 
  M_n(A)/(L+R)\cong (A^n/V(R))\otimes_A (A^n/V(L)). 
  $$
  \end{cor} 
 \proof
 By (\ref{V(L)_e}), (\ref{V(R)_e}) and (\ref{MN-decomposition}), respectively, we have
$$
\aligned
 M_n(A)&/(L+R)= \\ 
 &(A^n\otimes_AA^n)/(A^n\otimes_AV(L)+V(R)\otimes_AA^n)=\\ &\qquad \qquad\qquad\qquad\qquad\qquad\qquad(A^n/V(R))\otimes_A (A^n/V(L)). 
 \endaligned
 $$
  \qed

\subsection{Locally free modules}
Recall that an algebra is {\bf affine} if it is finitely generated as an algebra. An algebra is a {\bf domain} if it has no zero divisors. A module $M$ over the algebra $A$ is {\bf finitely generated} if  there is a finite subset $\{b_i\}_{i=1}^n$ of $M$ so that every element of $M$ can be written as an $A$-linear combination of the $\{b_i\}_{i=1}^n$. Finally, the module $M$ is free if the generating set $\{b_i\}_{i=1}^n$ can be chosen to be $A$-linearly independent. In this case, the cardinality of  $\{b_i\}_{i=1}^n$ is the {\bf rank}  of the free module $M$.  \\

In this section, we give a criterion for when certain modules over an affine domain can be localized so that they become free. The proof uses Nakayama's lemma and the fact that if $A$ is an affine commutative domain over the complex numbers then the intersection of all maximal ideals in $A$ is zero. \\

The following form of Nakayama's lemma appears in \cite[p.124]{A}.
\begin{theorem}[Azumaya] \label{bases}   Let $M$ be a finitely generated module over the commutative ring $A$ with identity. Let $N\leq M$ be a submodule of $M$ so that for every maximal ideal $\fm$ of $A$, $M=N+\fm M$, then $N=M$. \end{theorem}

 The algebra $K$ is {\bf prime} if given  $a,b\in K$, if for all $r\in K$,  $arb=0$ then $a=0$ or $b=0$. If $K$ is prime, then its center $Z(K)$ is a domain. If the algebra $K$ has no zero divisors then it is prime. 

\begin{prop} \label{free}Let $K$ be a finitely generated prime  algebra over the complex numbers.  Suppose that $A\leq Z(K)$ has the property that $K$ is a finitely generated module over $A$. Let $\fm\leq A$ be a maximal ideal with an open neighborhood $U$ in $\Spec(A)$ such that 
$$ 
 \forall \fm'\in U\cap \MS(A) \qquad \mathrm{dim}_\mathbb{C}(K/\fm'K)=d.
$$ 
Then there exist $\{b_i\}_{i=1}^d\in K$ and  $\alpha\in A$ with $\fm\in D(\alpha)$ such that  $K_\alpha$ is {\bf free} over $A_\alpha$ with a basis the image $\{\overline{b}_i\}_{i=1}^d$ of  $\{b_i\}_{i=1}^d$ in the localization. In other words,
$$ 
K_\alpha\cong\bigoplus_{i=1}^d A_\alpha \overline{b}_i.
$$ 
\end{prop}

\proof By the Artin-Tate lemma, $A$ is an affine commutative algebra over the complex numbers.  Since $K$ is prime, $A$ has no zero divisors.  Therefore, 
$$
\bigcap_{\fm \in \MS(A)}\fm=0.
$$

Since, by assumption, $\dim(K/\fm K)=d$, there is a subset $\{b_i\}_{i=1}^d\subset K$ whose image $\{\ov{b}_i\}_{i=1}^d\subset K/\fm K$ forms a basis for $K/\fm K$ as a vector space over $A/\fm=\mathbb{C}$. Being linearly independent is an open condition in the Zariski topology. 
Thus, there exists a Zariski open  set $V\subset \Spec(A)$  containing $\fm$ so that the image of $\{b_i\}_{i=1}^d$ is linearly independent in $K/\fm'K$ for all $\fm'\in V$; hence, the image of $\{b_i\}_{i=1}^d$ is a basis for $K/\fm'K$ for all $\fm'\in U\cap V$.\\

The principal open subsets form a basis for the Zariski topology. Thus,  there exists $\alpha \in A$ with $\fm \in D(\alpha)\subset V\cap U$.\\

By Theorem \ref{bases}, the image of $\{b_i\}_{i=1}^d$ in $K_\alpha$, denoted  $\{\overline{b}_i\}_{i=1}^d$,
spans $K_\alpha$. We need to check that $\{\overline{b}_i\}_{i=1}^d$ is linearly independent.  Suppose
\begin{equation}\label{badnews} 
\sum_ia_i\overline{b}_i=0,\quad \mathrm{with}~~a_i \in A_\alpha \quad \forall~i=1,\ldots,d,
\end{equation}
is a non-trivial relation. Since the intersection of all maximal ideals of $A_\alpha$ is zero, we conclude that if $a_i\neq 0$, there exists
$\fm'\in \MS(A_\alpha)$ so that the image of $a_i$ in $A/\fm'$ is nonzero. This implies that the image of $\{\overline{b}_i\}$ in $K/\fm'K$  is linearly dependent, which is a contradiction. Therefore, there is no non-trivial relation as in (\ref{badnews}) and $\{\overline{b}_i\}_{i=1}^d$ freely generates $K_\alpha$ as an $A_\alpha$ module.  \qed

\begin{rem}  Looking ahead to the application of this proposition,
note that, if $X$ and $\fm$ are as in Theorem~\ref{main_thm}, then we cannot use Proposition~\ref{free}  by applying it directly to the pair 
$$
(K,A)=(\Kz(X),\Ke(X))
$$ 
in order to prove  the Theorem~\ref{main_thm}, because  $\Ke(X)$ is (usually) not a domain.  In order to use this proposition  to derive the desired conclusion we need to  work with a Heegaard splitting of the $3$-manifold $X$, along with some facts about skein algebras of surfaces. 
\end{rem}

\subsection{The Azumaya Locus} Suppose that the algebra $K$  is finitely generated, prime, and is a finite rank module over its center. By the Artin-Tate lemma, its  center $Z(K)$  is a finitely generated integral domain. Let
$S$ be the multiplicatively closed subset of nonzero elements of $Z(K)$. The localization $S^{-1}Z(K)$ of the center of $K$ at $S$ is a field and $S^{-1}K$ is a finite dimensional algebra over that field. The {\bf dimension} of $K$ is the dimension of $S^{-1}K$ as a vector space over $S^{-1}Z(K)$.  Furthermore, by Posner's theorem \cite{MR}, the dimension  is $d^2$ for some natural number $d$.   The {\bf Azumaya locus} of $K$ is the set 
$$
\mathcal{A}(K)=\Big\{ \fm\in \MS(Z(K))\colon K/\fm \cong M_d(\mathbb{C})\Big\}\subset \MS(Z(A)).
$$
We say $K$ is {\bf Azumaya} if  $\mathcal{A}(K)=\MS(Z(K))$; c.f. \cite{MA}. Proposition \ref{free}, applied with $A=Z(K)$, shows that an Azumaya algebra is isomorphic to the global sections of a vector bundle over the variety of $Z(K)$ with fiber $M_d(\mathbb{C})$.  That is, $K$ is a projective module over $Z(K)$, \cite{Se}  \\

 \section{Towards the Kauffman bracket skein module of a closed three-manifold}\label{Skein_sec}
In this section, we continue to recall known results about the structure of the Kauffman bracket skein module of an oriented $3$-manifold, using the  terminology introduced in section \ref {Algebra_sec}. We conclude by reviewing the proof that the skein module $\Kz(M)$ for a $3$-manifold is a finite rank module over $\Ke(M)$.  \\

Let $X$ be an oriented closed $3$-manifold. Two representations 
 $
 \rho_1,\rho_2\colon \pi_1(X)\rightarrow \G
 $ 
 are {\bf trace equivalent} if
$$ 
Tr(\rho_1(\gamma))=Tr(\rho_2(\gamma))\qquad \forall \gamma \in \pi_1(X).
$$
The skein algebra $\Ke(X)$ is affine, and its maximal spectrum
   is in one-to-one correspondence with trace equivalence classes of representations 
 $$
 \rho\colon \pi_1(X)\rightarrow \G
 $$ 
(see Theorem \ref{Dougsbigthm}).\\

Let $\Sigma_{g}$  be a closed oriented surface of genus $g$. For $\zeta$ and $\epsilon$ as specified in (\ref{mep_e}), the threading map $\tau\colon \Ke(\Sigma_{g})\rightarrow Z(\Kz(\Sigma_{g}))$ in (\ref{tau}) is an isomorphism and allows us to identify $\MS(Z(\Kz(\Sigma_{g})))$ with $\MS(\Ke(\Sigma_{g}))$.

\begin{theorem}\cite{GJS2,KK}. 
With the notation as above, the Azumaya locus $\mathcal{A}(\Kz(\Sigma_{g}))$ consists of maximal ideals in $\MS(\Ke(\Sigma_{g}))$ that do not correspond to central representations.
\end{theorem}
The Azumaya locus of the skein algebra of a surface is an open set. Note that for every $\fm\in \mathcal{A}(\Ke(\Sigma_{g}))$ we have 
 $$
 \Kz(\Sigma_{g})/\tau(\fm)\Kz(\Sigma_{g})\cong M_d(\mathbb{C})
 $$  
 where $d=m^{3g-3}$ \cite{FKL1}.\\

Moving up in dimension from surfaces to $3$-manifolds, suppose $H$ is a handlebody with $\partial H\!=\!\Sigma_{g}$. It is possible to find a properly embedded planar surface
   $S\!\subset\! H$ so that $H\cong S\times [0,1]$. Such surface $S$ is called a {\bf spine} of the handlebody.  Note that 
   $$
   \Sigma_{g}= S\times \{0,1\}\cup \partial S\times [0,1].
   $$
   Figure \ref{h2} illustrates a genus two handlebody as a cylinder over a planar surface $\Sigma_{0,3}$.
   \begin{figure}[H]
   \includegraphics[scale=.1]{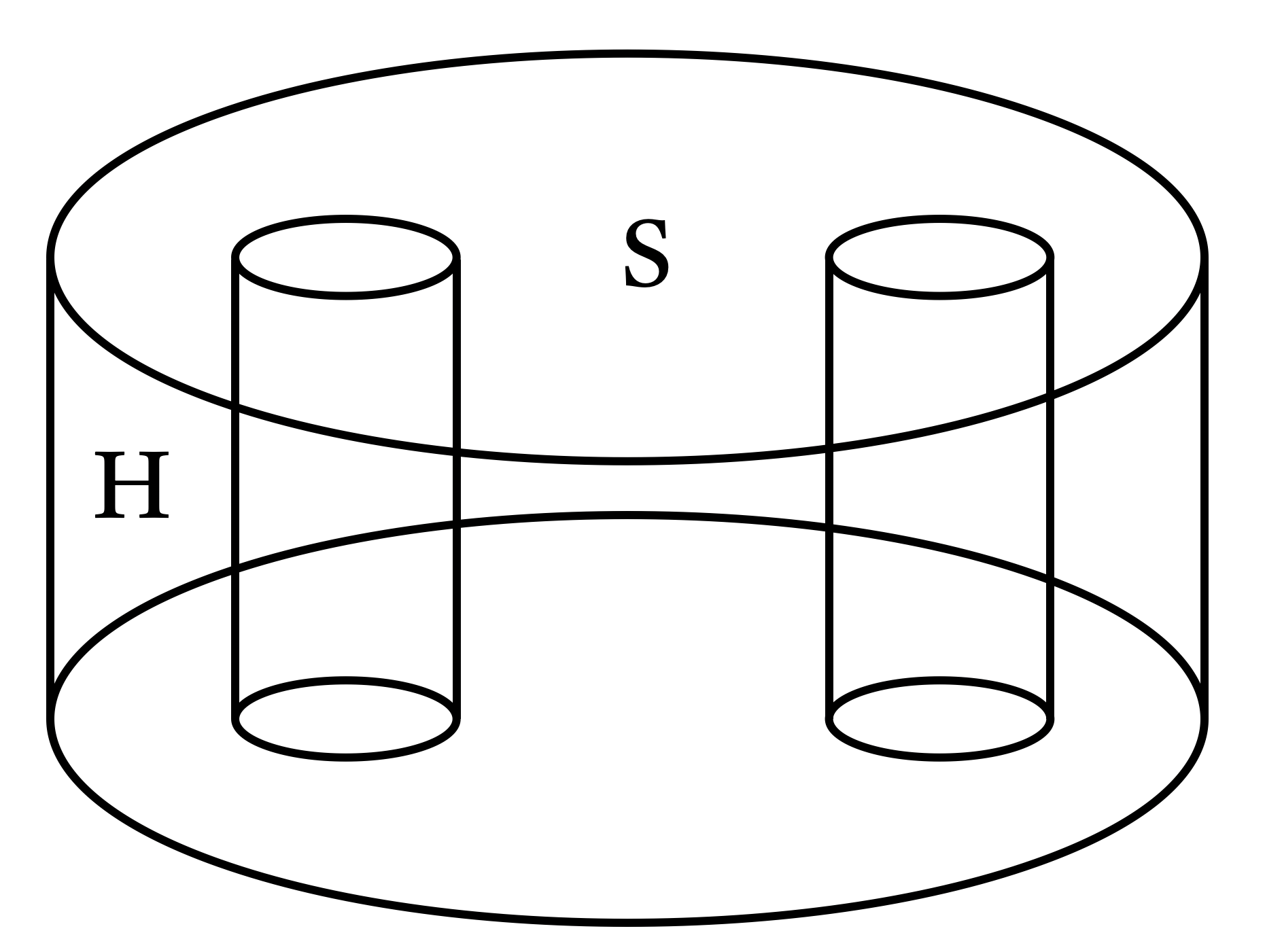} \caption{Handlebody of genus $2$}\label{h2}
   \end{figure} 
A planar spine $S$ for a genus $g$ handlebody $H$ has $g+1$ boundary components; i.e. $S=\Sigma_{0,g+1}$. Any framed link in $H$ can be pushed into a collar of the boundary lying over the spine  $\Sigma_{0,g+1}  \times \{1\}$, or pushed  into a collar of the boundary lying under the spine $\Sigma_{0,g+1}\times \{0\}$.  Depending on whether you imagine $\Kz(H)$ to be a left or a right module, one of these pushes yields an algebra morphism
\be
\sigma_S:\Kz(\Sigma_{0,g+1} )=\Kz(H)\rightarrow \Kz(\Sigma_{g}).
\ee
In the other direction, the map
\be
h\colon \Kz(\Sigma_{g})\rightarrow \Kz(H)
\ee
comes from identifying $\Sigma_{g} \times [0,1]$ with a collar of $\Sigma_{g}=\partial H$ in $H$ (also in one of two possible ways).
A composition of these two  maps
satisfies
\be
h\circ \sigma_S=\mathrm{id}_{\Kz(H)}. 
\ee
There are right and left ideals  $R_H$  and  $L_H$  in $\Kz(\Sigma_{g})$, respectively, depending on how $\Sigma_{g} \times [0,1]$ is identified with a collar of $\partial H$ in $H$, that are the annihilator of the empty skein in $H$ and fit into the short exact  sequences of $\Ke(\Sigma_{0,g+1})$-modules.  Here, $h_L$ and $h_R$ are the corresponding choices of $h$. 

\be  \label{SESofR} 
\begin{tikzcd} 0\arrow{r}& R_H\arrow{r}{\iota_R} & \Kz(\Sigma_{g}) \arrow{r}{h_R}  & \Kz(H) \arrow[bend left=33]{l}{\sigma_S}  \arrow{r} & 0 \end{tikzcd} 
\ee
 or 
  \be \label{SESofL}
   \begin{tikzcd} 0\arrow{r}& L_H\arrow{r}{\iota_L} & \Kz(\Sigma_{g}) \arrow{r}{h_L}  & \Kz(H) \arrow[bend left=33]{l}{\sigma_S}  \arrow{r} & 0. \end{tikzcd}
 \ee
Note, however, that  viewed as  sequences of $\Ke(\Sigma_g)$-modules, the sequences above do not split since  the skein modules of the subsurfaces lying above and below $\Sigma_{0,g+1}$ in $\Sigma_g$ are not stabilized by the action of $\Ke(\Sigma_g)$. \\

The Azumaya locus of $\Kz(\Sigma_{0,b})$ contains all maximal ideals corresponding to irreducible representations of $\pi_1(\Sigma_{0,b})$ such that none of the peripheral curves gets sent to matrices of trace $\pm 2$; c.f. \cite{FKL}. We need a little bit more.

 \begin{theorem}[\cite{FKL}]\label{bigfish} Let $H$ be a handlebody of genus $g$.
 Suppose $\rho\colon\pi_1(H)\rightarrow \G$ is irreducible and $\fm \in \MS(\Ke(H))$  is the maximal ideal corresponding to $\rho$. 
 There exists a spine $\Sigma_{0,g+1}\subset H$ such that the pullback of $\rho$ to $\pi_1(\Sigma_{0,g+1})$ has the property that no peripheral loop is sent to a matrix of trace equal to $\pm 2$. Furthermore,
   $$ 
   \Kz(H)/\tau(\fm)\Kz(H) 
   $$ 
   is a vector space of dimension $m^{3g-3}$ and there is a Zariski open neighborhood  $U$ of  $\fm$ in $\MS(\Ke(H))$ so that for every $\fm'\in U$, 
   $$ 
   \Kz(H)/\tau(\fm')\Kz(H) 
   $$ 
   also has dimension $m^{3g-3} $. \end{theorem}

 \proof The only point not addressed in \cite{FKL} is that there is an open neighborhood of $\fm$ so that the quotients are complex vector spaces of dimension $m^{3g-3}$. The conditions that the representation are irreducible, and that none of the peripheral loops are sent to matrices of trace $\pm 2$ are open. \qed\\

Given a closed oriented $3$-manifold $X$, suppose   $H,B\subset X$  are handlebodies that intersect along their joint boundary $F$,  and  $X=H\cup B$. We say $(H,B,F)$ is a  {\bf  Heegaard splitting} of $X$. Every closed oriented three-manifold admits a Heegaard splitting. If $(H,B,F)$ is  a Heegaard splitting of $X$ then 
  \be \label{tensorr} 
  \Kz(X)=\Kz(H)\otimes_{\Kz(F)}\Kz(B); \ee 
  c.f. \cite{Przy}. We highlight the following statement from  \cite{FKL} as it is central to this paper.

  \begin{theorem}[\cite{FKL}] With $\zeta$ and $\epsilon$ as in (\ref{mep_e}),  for every  closed three-manifold $X$,  $\Kz(X)$ is a finitely generated module over $\Ke(X)$. \end{theorem}

  \proof Fix a Heegaard splitting $(H,B,F)$.  By (\ref{tensorr}), (\ref{SESofR})-(\ref{SESofL}), and Corollary~\ref{L+R_cor}, we have
  \be \label{Kzeta/L+R} 
  \Kz(X)=\Kz(F)/(L_B+R_H)
   \ee
where $R_H$ and $L_B$ are as in  (\ref{SESofR})-(\ref{SESofL}).
By \cite{AF}, $\Kz(F)$ is a finitely generated module over $\Ke(F)$. Consequently, the quotient (\ref{Kzeta/L+R}) is also a finitely generated module over $\Ke(F)$. By the naturality of the threading map, the action of $\Ke(F)$ on the quotient factors through $\Ke(X)$; thus, $\Kz(X)$ is a finite rank module over $\Ke(X)$.  \qed

\section{Localized skein module modules at an irreducible representation}\label{localized-model_sec} 
In this section we state the main theorem of the paper and lay  groundwork for the proof, which is given in the following section \ref{proof_sec}.\\

In a recent paper \cite{FKL} it was proved that if $X$ is a closed oriented homology sphere, $\rho:\pi_1(X)\rightarrow \G$ is an irreducible representation,  and $\fm\leq \Ke(X)$ is the maximal ideal corresponding to $\rho$, then the specialization  $\Kz(X)/\fm\Kz(X)$ is a vector space of dimension one over $\Ke(X)/\fm\Ke(X)\cong \mathbb{C}$. The main result of this paper is the following enhancement of the aforementioned  result.\vskip.1in

\begin{theorem}\label{main_thm} Suppose that $\zeta$ is a primitive $n$-th root of unity with  $n$ not divisible by $4$,  and let $\epsilon=\zeta^{(\frac{n}{ \gcd(n,2)})^2}.$ 
Let $X$ be a closed oriented $3$-manifold, and  $\fm\leq \Ke(X)$  the maximal ideal corresponding to an irreducible representation $\rho:\pi_1(X)\rightarrow \G$. The localized skein module $\Kz(X)_\fm $  is a free module of rank one over $\Ke(X)_\fm$; i.e. $\Kz(X)_\fm\cong \Ke(X)_\fm$. 
\end{theorem}

The proof proceeds along the following outline. 
Consider a Heegaard splitting $(H,B,F)$  of  $X$ obtained by cutting $X$ into two handlebodies $H$ and $B$ along a genus $g$  surface $F$.  The skein modules $\Kz(H)$ and $\Kz(B)$ are  right and left modules, respectively, over the skein algebra $\Kz(F)$ of the surface $F$. 
Let $d=m^{3g-3}$ where $m= \frac{n}{ \gcd(n,2)}$ as in (\ref{mep_e}). 
\begin{itemize}[wide, nosep, labelindent = 0pt, topsep = 1ex]
\item In \cite{GJS2}, it is proved that $\Kz(F)$ specialized at the maximal ideal $\fm$ of an irreducible representation $\rho$ is isomorphic to the complex vector space $M_d(\mathbb{C})$ of $d\times d$ matrices;\vskip.1in
\item In \cite{FKL} it is proved that $\Kz(H)$ and $\Kz(B)$ specialize at $\fm$ to $\mathbb{C}^d$. 
\end{itemize}

Consequently, $\Kz(X)$ specialized at $\fm$
is isomorphic to
\begin{multline} 
  \Kz(X)/\tau(\fm)\Kz(X) \cong\\
  \left(\Kz(H)/\tau(\fm)\Kz(H)\right)\otimes_{\Kz(F)/\tau(\fm)\Kz(F)}\left(\Kz(B)/\tau(\fm)\Kz(B)\right) \cong \\ \mathbb{C}^d\otimes_{M_d(\mathbb{C})}\mathbb{C}^d =\mathbb{C},
\end{multline} 
see \cite[Theorem 12.2]{FKL}. The latter follows from the fact that the tensor product of the row vector representation and the column vector representation
of a matrix algebra is isomorphic to $\mathbb{C}$.\\

The conditions that the specializations of the algebras and skein modules have dimension $d^2$ and $d$, respectively, are open in the Zariski topology. Using Nakayama's lemma we are able to prove that localizations in small enough neighborhoods of the skein algebra and the skein modules are free modules of the appropriate dimensions and then carry out the analogous argument to prove our theorem.

\begin{rem}Our technique of proof fails at the trivial representation.
 \end{rem}

Fix an irreducible representation $\rho:\pi_1(X)\rightarrow \G$ and denote the corresponding maximal ideal of $\Ke(X)$ by $\fm$. Also choose  a Heegaard splitting $(H,B,F)$ of $X$.
   The inclusions $j\colon B\rightarrow X$ and $i\colon H\rightarrow X$, together with $h_R$ associated to $H$ as in (\ref{SESofR}) and $b_L$ associated to $B$ as in (\ref{SESofL}), induce surjective algebra morphisms fitting in a commutative diagram
$$ 
\begin{CD}  \Ke(F) @>h_R>> \Ke(H) \\ @V b_L VV @ViVV \\ \Ke(B) @>j>> \Ke(X). \end{CD} 
$$ \vskip.1in
By the surjectivity of these maps, we can identify  
$$
\Spec(\Ke(X)), \quad \Spec(\Ke(H)),\quad\mathrm{and}\quad \Spec(\Ke(B))
$$ 
with subspaces of $\Spec(\Ke(F))$.  Therefore, the maximal ideals corresponding to $\rho$ and its pullbacks to $\pi_1(H)$, $\pi_1(B)$ and $\pi_1(F)$  will all be denoted by the same notation~$\fm$. 
 %The assumption that $\rho$ is an isolated  point can be rephrased as:   $\fm$ is an isolated intersection point of the holomorphic Lagrangians $\Spec(\Ke(H))$ and $\Spec(\Ke(B))$ in the holomorphic symplectic variety $\Spec(\Ke(F))$. Furthermore, $\dim_{\mathbb{C}} \Ke(F)_\fm$ is the order of contact of the two Lagrangians at $\fm$ (we will not make use of these observations in the following arguments).  

 \begin{rem}  \label{localization-notation_rmk}
 The homomorphisms $h_R$, $b_L$, and $j\circ b_L=i\circ h_R$ make $\Ke(H)$, $\Kz(H)$, $\Ke(B)$, $\Kz(B)$, $\Ke(X)$, and $\Kz(X)$  into modules over $\Ke(F)$.  For instance, localizing $\Ke(H)$ as a module over $\Ke(F)$ at $\alpha$ has the same effect as localizing $\Ke(H)$ as an algebra at $h_R(\alpha)$. That is,
 $$
 \Ke(H)_{h_R(\alpha)}=\Ke(H)_\alpha. 
 $$ 
 The notation on the right is cleaner.  When we are localizing at the powers of an element  of these algebras, we denote it by sub-scripting with  a preimage  of the element in $\Ke(F)$.\end{rem}

\begin{prop} \label{freef} Given a Heegaard splitting $(H, B, F)$ of a $3$-manifold $X$ and  with notation as above,
\begin{enumerate}[wide, nosep, labelindent = 0pt, topsep = 1ex]
    \item there is $\delta\in \Ke(F)$ with $\fm\in D(\delta)$ so that $\Kz(F)_\delta$ is a free module over $\Ke(F)_\delta $ of rank $d^2=m^{6g-6}$; implying $\Kz(F)_\fm$ is a free module over $\Ke(F)_\fm$ of rank $d^2$.
    \item there is $\delta\in \Ke(F)$ with $\fm \in D(\delta)$ so that $\Kz(H)_\delta$ is a free module over $\Ke(H)_\delta$ of rank $d$; implying that $\Ke(H)_\fm$ is a free module of rank $d$. 
      \end{enumerate} 

   \end{prop} 

   \proof  Since $\rho$  is irreducible, the ideal $\fm$  is in the Azumaya locus of $\Ke(F)$.  The Azumaya locus of $\Ke(F)$ is open, hence there  is an open subset $U_F'\subset \Spec(\Ke(F))$ such that for all maximal ideals $\fm' \in U_F'$, 
   $$ 
   \Kz(F)/\tau(\fm')\Kz(F) \cong M_d(\mathbb{C}),
   $$ 
   where $d=m^{3g-3}$. Since $\tau:\Ke(F)\rightarrow \Kz(F)$ is injective, and the algebra $\Kz(F)$ has no zero-divisors and  is a finitely generated module over its center, Proposition \ref{free} applies where $\Kz(F)$ plays the role of the algebra $K$. The subalgebra $A$  is the center of $\Kz(F)$.  Therefore  there is $\delta \in \Ke(F) $, with  $\fm \in D(\delta)$, so that  $\Kz(F)_\delta$ is a free module over $\Ke(F)_\delta$ of rank $d^2$.\\

Recall the discussion of the limit of a directed sum of modules, with the
 map $\Ke(F)_\delta\rightarrow \Ke(F)_\gamma$  as in Equation \ref{directed}, where  the neighborhood basis  $D(\gamma)\subset D(\delta)$ is used to make the localizations of $\Ke$ into a directed system.  We claim any basis  $\{b_i\}$ for $\Kz(F)_\delta$ as a free module over $\Ke(F)_\delta$,  is sent
   to a basis for $\Kz(F)_\gamma $ as a module over $\Ke(F)_\gamma$ by the map
   from the directed system because it can be understood as extension of scalars. Notice that if there is a  nontrivial linear dependence among the image of the $\{b_i\}$ in $\Kz(F)_\fm$ we can clear the fractions to get a nontrivial dependence between the $\{b_i\}$ in $\Kz(F)_\delta$, which is a contradiction. \\

   Letting $I=\{D(\gamma)\}_{\fm \in D(\gamma)}$,   analogously to (\ref{dirmod}), we have 
   \be \Kz(F)_\fm=\varinjlim_{D(\gamma)\in I} \Kz(F)_\gamma. \ee
 From the definition of the direct limit of a directed system of modules, the image of $\{b_i\}$ in $\Ke(F)_\fm$ spans.   Once again, given a  linear dependence among the image of the $\{b_i\}$ in $\Kz(F)_\fm$, using the fact that $\Ke(F)$ has no zero divisors we can clear the fractions in the linear dependence to get a linear dependence among the $\{b_i\}$ in $\Kz(F)_\gamma$ for some $\gamma$. This contradicts the fact that the image of the $\{b_i\}$ in $\Ke(F)_\gamma$ forms a basis.  \\
  
To prove the second part, choose a planar spine $\Sigma_{0,g+1}\subset H$ where $g$ is the genus of $F$ so that none of the peripheral curves is sent to a matrix with trace $\pm 2$ by the representation corresponding to $\rho$. By Theorem \ref{bigfish}, there is a Zariski open neighborhood $U$ of $\fm$ in $\MS(\Ke(H))$ so that for all $\fm'\in U$,
  $$ 
  \dim_{\mathbb{C}}\Kz(H)/\fm'\Kz(H)= d.
  $$  
  We now apply Proposition \ref{free} using $\Kz(\Sigma_{0,g+1})$ as the algebra. The role of the subalgebra of the center of $\Kz(\Sigma_{0,g+1})$ is played by $\tau(\Ke(\Sigma_{0,g+1}))$.  Since 
  $$
  \Spec(\Ke(H))\subset\Spec(\Ke(F))
  $$ 
  is an embedding, we can choose $\delta\in \Ke(F)$ so that $\fm \in D(\delta)$ and $\Kz(H)_\delta$ is a free module of rank $d$ over $\Ke(H)_\delta$. We take the limit similarly to above to see that $\Kz(H)_\fm$ is a free module of rank $d$ over $\Ke(H)_\fm$.
 \qed

 \begin{prop}\label{matrix} Suppose that $F=\partial H$ where $H$ is a handlebody of genus $g$. Suppose that
   $\fm\in \MS(\Ke(F))$ corresponds to the pullback of an  irreducible representation of $\pi_1(H)$.  There exists
   $\delta \in \Ke(F)$ with $\fm \in D(\delta)$ such that
   $$ 
   (\Kz(F)\otimes_{\Ke(F)}\Ke(H))_\delta \cong M_d(\Ke(H)_\delta). 
   $$ This implies that $$ \left(\Kz(F)\otimes_{\Ke(F)}\Ke(H)_\fm\right)\cong M_d(\Ke(H)_\fm).$$
   \end{prop}

 \proof 
 Let $\fp \in\Spec(\Ke(F))$ be the kernel of the homomorphism $h:\Ke(F)\rightarrow \Ke(H)$ coming from inclusion. The ideal $\fp$ is prime because the character variety of a free group is irreducible.  Note that
 $$ 
 \Kz(F)\otimes_{\Ke(F)}\Ke(H)=\Kz(F)/\tau(\fp)\Kz(F). 
 $$
 This means that the  (left or right) action of $\Kz(F)$ on $\Kz(H)$ descends to give an action
 $$ 
 \left(\Kz(F)\otimes_{\Ke(F)}\Ke(H)\right) \otimes \Kz(H)\rightarrow \Kz(H). 
 $$
 This in turn gives a homomorphism
 $$ 
 \theta:\Kz(F)\otimes_{\Ke(F)}\Ke(H) \rightarrow \End_{\Ke(H)}(\Kz(H)). 
 $$
 Specializing at the point $\fm$, this becomes an isomorphism 
 $$
 \Kz(F)/\fm\Kz(F) \stackrel{\cong}{\longrightarrow}\End_{\mathbb{C}}\big(\Kz(H)/\fm\Kz(H)\big).
 $$ 

By Proposition \ref{freef}, there is  $\delta_1\in \Ke(F)$ so that $\fm \in D(\delta_1)$ and $\Kz(F)_{\delta_1}$ is a free $\Ke(F)_{\delta_1}$ module of rank $d^2$. Similarly choose $\delta_2 \in \Ke(F)$ with $\fm \in D(\delta_2)$ so that $\Kz(H)_{\delta_2}$ is a free module of rank $d$ over $\Ke(H)_{\delta_2}$.\\

 Choosing $\{b_i\}_{i=1}^{d^2}$ in $\Kz(F)$ that form a basis of
 $\Kz(F)/\fm\Kz(F)$, the image of this basis is   linearly independent in an open neighborhood of $\fm$.
Therefore, there is $\delta$ with $\fm \in D(\delta)\subset D(\delta_1)\cap D(\delta_2)$ so that
 $$ 
 \theta:\left(\Kz(F)\otimes_{\Ke(F)}\Ke(H)\right)_\delta \rightarrow \End_{\Ke(H)_\delta}\Kz(H)_\delta 
 $$ 
 is an isomorphism. However, after choosing a basis for the free module $\Kz(H)_\delta$ over $\Ke(H)_\delta$, we have
 $$
 \End_{\Ke(H)_\delta}\Kz(H)_\delta=M_d(\Ke(H)_\delta).
 $$

 Since $\theta$ remains an isomorphism over sufficiently small principal open neighborhoods, we get that
 $$(\Kz(F)\otimes_{\Ke(F)}\Ke(H))_\fm\cong M_d(\Ke(H)_\fm). $$ \qed

 \section{Proof of Theorem~\ref{main_thm}}\label{proof_sec}
We use the results established in Section \ref{localized-model_sec} to conclude the proof of (the main) Theorem \ref{main_thm}.
\proof 
Recall that $X$ is a closed oriented $3$-manifold. Starting from a Heegaard splitting $(H,B,F)$ of $X$, by (\ref{tensorr}) and (\ref{Kzeta/L+R}), we have
    $$ 
    \Kz(X)=\Kz(H) \otimes_{\Kz(F)} \Kz(B)= \Kz(F)/(R_H+L_B). 
    $$
    If $\fp_H,\fp_B\leq \Ke(F)$ are the prime ideals that are the kernels of the maps induced by inclusions
    $$ h_R:\Ke(F)\rightarrow \Ke(H) \quad \mathrm{and} \quad b_L:\Ke(F)\rightarrow \Ke(B),
    $$
    respectively, then
    $$ 
    \fp_H\Kz(F)\leq R_H \quad \mathrm{and} \quad \fp_B\Kz(F)\leq L_B. 
    $$
Dividing both the numerator and the denominator of the right hand side of (\ref{Kzeta/L+R}) by $(\tau(\fp_H)+\tau(\fp_B))\Kz(F)$, we get
\be\label{bigfrac}  
\Kz(X)=\left(\frac{\Kz(F)}{(\tau(\fp_H)+\tau(\fp_B))\Kz(F)}\right)/\left(\frac{R_H+L_B}{(\tau(\fp_H)+\tau(\fp_B))\Kz(F)}\right).
\ee
    Taking the quotient by $(\tau(\fp_H)+\tau(\fp_B))\Kz(F)$ is the same as tensoring with 
    $$
    \Ke(F)/(\fp_H+\fp_B)=\Ke(X)
    $$ 
    over $\Ke(F)$. Therefore,
   \be\label{top}
   \Kz(F)/(\tau(\fp_H)+\tau(\fp_B))\Kz(F)=\Kz(F)\otimes_{\Ke(F)}\Ke(X)
   \ee 
   and
  \be\label{bottom} \frac{(R_H+L_B)}{(\tau(\fp_H)+\tau(\fp_B))\Kz(F)}=R_H\otimes_{\Ke(F)}\Ke(X)+L_B\otimes_{\Ke(F)}\Ke(X).\\
  \ee
  Substituting \ref{top} and \ref{bottom} in \ref{bigfrac}, this means 
  \be \label{almost} \Kz(X)=\left(\Kz(F)\otimes_{\Ke(F)} \Ke(X)\right)/\left(R_H\otimes_{\Ke(F)}\Ke(X) +L_B\otimes_{\Ke(F)}\Ke(X)\right). \ee

Since tensor product is a right exact functor, the result of tensoring the  short exact sequences  (\ref{SESofR}) associated to $H$ and (\ref{SESofL}) associated to $B$ with $\Ke(H)$ and $\Ke(B)$, yields  exact sequences of $\Ke(H)$-modules and $\Ke(B)$-modules,
 \be \label{forH} \begin{CD}
 R_H\otimes_{\Ke(F)}\Ke(H)@>\iota_R>> \Kz(F)\otimes_{\Ke(F)}\Ke(H)@>h_R>> \Kz(H)\rightarrow 0 \end{CD}
 \ee 
 and 
 \be \label{forB}  \begin{CD}
  L_B\otimes_{\Ke(F)}\Ke(B)@>\iota_L>> \Kz(F)\otimes_{\Ke(F)}\Ke(B)@>b_L>> \Kz(B)\rightarrow 0\end{CD},
 \ee
 respectively.
 Here, on the right end of each exact sequence,  we are using the fact that if you tensor a module by its base ring you recover the module.\\ 

 Next, by Propositions \ref{freef} and \ref{matrix} and the fact that the principal open sets form a basis for the Zariski topology, there exists $\delta\in \Ke(F)$ so that
 $$ \Kz(F)_\delta\cong M_d(\Ke(F)),
 $$
 $$ 
 \Kz(H)_\delta\cong \Ke(H)^d \quad \mathrm{and}\quad \Kz(B)_\delta=\Ke(B)^d .
 $$ 
 Recall from (\ref{tensorr}) that
 $$ 
 \Ke(X)=\Ke(H)\otimes_{\Ke(F)}\Ke(B).
 $$ 
 Tensor product commutes with localization. Tensoring Equation \ref{forH} by $\Ke(B)$ and Equation (\ref{forB}) by $\Ke(H)$, followed by localizing at the principle ideal generated by $\delta$, we get
 \be \label{forH1} 
 R_H\otimes_{\Ke(F)}\Ke(X)_\delta \rightarrow M_d(\Ke(X)_\delta) \rightarrow \Ke(X)^d_\delta\rightarrow 0 
 \ee 
 and 
 \be \label{forB1}  
  L_B\otimes_{\Ke(F)}\Ke(X)_\delta\rightarrow M_d(\Ke(X)_\delta)\rightarrow \Ke(X)_\delta^d \rightarrow 0. 
 \ee

 Taking  limits  we find that
 \be \label{forH2} 
  R_H\otimes_{\Ke(F)}\Ke(X)_\fm \rightarrow M_d(\Ke(X)_\fm) \rightarrow \Ke(X)^d_\fm\rightarrow 0 
 \ee 
 and 
 \be \label{forB2}  
  L_B\otimes_{\Ke(F)}\Ke(X)_\fm\rightarrow M_d(\Ke(X)_\fm)\rightarrow \Ke(X)_\fm^d \rightarrow 0. 
 \ee

 If $V(R)$ is the row space of the image of  $ R_H\otimes_{\Ke(F)}\Ke(X)_\fm$ in $ M_d(\Ke(X)_\fm)$  and $V(L)$ is the column space of the image of 
 $ L_B\otimes_{\Ke(F)}\Ke(X)_\fm$ in $  M_d(\Ke(X)_\fm)$  we conclude from Equations \ref{forH2} and \ref{forB2} that
 \be \label{sothere} ((\Ke(X)_\fm)^d/V(R))^d\cong \Ke(X)^d \ \mathrm{and} \  ((\Ke(X)_\fm)^d/V(L))^d\cong \Ke(X)^d.\ee
 This implies that $\Ke(X)_\fm^d/V(R)$ and $\Ke(X)_\fm^d/V(L)$ are projective modules over $\Ke(X)_\fm$.  However, by \cite{Ka}, any projective module over a local ring is free. Since the direct sum of $d$ copies of the free modules  $((\Ke(X)_\fm)^d/V(R)$ and $(\Ke(X)_\fm)d/V(L)$ are isomorphic to the direct sum of $d$ copies of the free module  $\Ke(X)_\fm$, we conclude that as $\Ke(X)_\fm$-modules
 $$ 
 (\Ke(X)_\fm)^d/V(R))\cong \Ke(X)_\fm \quad \mathrm{and} \quad  (\Ke(X)_\fm)^d/V(L))\cong \Ke(X)_\fm.
 $$ 
 Localizing Equation \ref{almost}  yields
 \begin{multline} \Kz(X)_\fm=\left(\Kz(F)\otimes_{\Ke(F)} \Ke(X)\right)/\left(R_H\otimes_{\Ke(F)}\Ke(X) +L_B\otimes_{\Ke(F)}\Ke(X)\right)_\fm =\\
   \Ke(X)^d_\fm/V(R)\otimes_{\Ke(X)_\fm}\Ke(X)^d_\fm/V(L) =\Ke(X)_\fm\otimes_{\Ke(X)_\fm}\Ke(X)_\fm =\Ke(X)_\fm.
\end{multline}\qed
 
\section{Applications and further directions}\label{cons} 
We use our result to improve an estimate of \cite{DKS} of the dimension of $\Kz(X)$ as a complex vector space. We also  estimate the dimension of  of $K_q(X)$ as a $\mathbb{Q}(q)$-vector space.
Throughout this section, assume that the $3$-manifold $X$ is a homology sphere with finite character variety.

\subsection{Estimating the dimension of $\Kz(X)$}

 When the order of the root of unity $\zeta$ is not divisible by $4$, 
we prove that the  dimension of the vector space $\Kz(X)$ is greater than or equal  to the dimension of the unreduced coordinate ring of the character variety as a complex vector space. This leads to a computation of the dimension of the skein module $K_{q}(X)$ with coefficients in rational functions, assuming that the universal skein module is tame.\\

Our goal here is to show how exploiting the fact that $\Kz(X)$ is a module over $\Ke(X)$ can lead to a deeper understanding of skein modules. The application  is based on arguments of \cite{DKS}.  As noted by the authors of that paper, a stronger result  can be proven by extending their arguments, i.e., if $X$ is a rational homology sphere with tame universal skein module then
$$
\dim_{\mathbb{Q}(q)} K_q(X)=\dim_\mathbb{C} \Ke(X). 
$$\\

The fact that $X$ is a homology sphere with finite character variety implies that
\be \label{natural-artinian-decomposition_e}
\Kz(X)=\Kz(X)_{\fm_0}\times \prod_{i=1}^n \Kz(X)_{\fm_i} 
\ee
where $\fm_0\in \MS(\Ke(X))$ is the maximal ideal corresponding to the trivial representation and $\{\fm_i\}_{i=1}^n$ are the maximal ideals corresponding to the irreducible representations.
Denote the quotient maps by 
\be\label{quot}
 q_i: \Kz(X)\rightarrow \Kz(X)_{\fm_i}.
\ee

 The Reshetkihin-Turaev invariant 
 $$ 
 RT^\zeta:\Kz(X)\rightarrow \Kz(S^3)=\mathbb{C} 
 $$ 
 is defined for $\zeta$ having order $2m$ where $m$ is odd. Let $t_{0}:\Ke(X)\rightarrow \mathbb{C}$ denote the algebra morphism having $\fm_0$ as its kernel.  If $L\leq X$ is a framed link with $c$ components, then for the skein represented by $L$ we have
    $$
     t_{0}(L)=(-2)^c. 
     $$
 The morphism $t_0$ makes $\mathbb{C}$ into a $\Ke(X)$-module.

 \begin{prop}[\cite{DKS} Proposition 2.4] The linear functional $RT^\zeta:\Kz(X)\rightarrow \mathbb{C}$ is a surjective intertwiner of $\Ke(X)$-modules where the module structure on $\mathbb{C}$ comes from~$t_{0}$. \qed\end{prop}

\begin{prop}\label{allgood}
If a primitive root of unity $\zeta$ has odd order $n$, then $-\zeta$ has order $2n$ and 
$$
\dim_{\mathbb C} \Kz(X)=\dim_{\mathbb C}K_{-\zeta}(X) 
$$
\end{prop}

 \proof It follows from a  theorem of John Barrett \cite{Ba}
 that every spin structure $s$ on $X$ induces a $\mathbb{C}$-linear  isomorphism
 $$ 
 \phi_s:\Kz(X)\rightarrow K_{-\zeta}(X).
 $$ \qed

We  prove the following.
  \begin{theorem}\label{itsbig} Let $X$ be a homology $3$-sphere with finite character variety.     
   Let $\zeta$ be a primitive root of unity of order $2m$, where $m$ is odd.   Let $q_i$ be the quotient maps defined by (\ref{quot}).
   The map 
   $$ 
   (RT^\zeta,q_1,\ldots,q_n):\Kz(X)\rightarrow \mathbb{C}\times \prod_{i=1}^n\Kz(X)_{\fm_i} 
   $$    
   is onto.  Therefore 
   $$
   \dim _{\mathbb C}\Kz(X)\geq  \dim _{\mathbb C}\Ke(X).
   $$
    \end{theorem}

    \proof The proof is an elaboration of the proof of Theorem 2.1 of \cite{DKS}. Recall   that 
$$
\Ke(X)\cong \Ke(X)_{\fm_0}\times \prod_{i=1}^n\Ke(X)_{\fm_i}
$$  
where $\fm_0$ is the maximal ideal corresponding to  the trivial representation, and for $i\geq 1$, the $\fm_i$ are the maximal ideals corresponding to  the irreducible representations. The map
\be (q_1,\ldots,q_n):\Kz(X)\rightarrow \prod_{i=1}^n\Kz(X)_{\fm_i} \ee is onto by the Chinese remainder theorem.
To complete the proof, we need to find an element of $\Kz(X)$ that is sent to zero by $ (q_1,\ldots,q_n)$ with nonzero image under $RT^\zeta$. 

  Let $t_i:\Ke(X)\rightarrow \mathbb{C}$ be the algebra morphism with kernel $\fm_i$. For every $i\neq 0$,
  since the representation corresponding to $\fm_i$ is irreducible, there exists a knot $K_i\subset X$ such that the skein of $K_i$ in $\Ke(X)$ has  $t_i(K_i) \not\in\{\pm 2\}$; see \cite{FKL}. Since $\Ke(X)$ is Artinian, for any maximal ideal  $\fm_i$ of $\Ke(X)$, 
there exists $k_i\in \mathbb{N}$ so that  if $c\in \fm_i$, then $q_i(c^{k_i})=0$. Let $k=\max \{k_i\}$. Consider the skein 
$$
 \varrho = \Big(\prod_{i=1}^n\big(K_i-t_i(K_i)\big)\Big)^k\in \Ke(X).
 $$

  Since   $\tau(\varrho)\in\fm_i^k$ for all $i>0$, we have   $q_i(\tau(\varrho))=0$.  Also $t_0(\varrho)\neq 0$. Choose a skein $\alpha\in \Kz(X)$ so that $RT^\zeta(\alpha)=1$. The skein
  $$ 
  \tau(\varrho)*\alpha  \in \Kz(X),
  $$
   where the $*$ denotes disjoint union and $\tau:\Ke(X)\rightarrow \Kz(X)$ is the threading map,
  satisfies
  $$ q_i(  \tau(\varrho)*\alpha)=0\in \Kz(X)_{\fm_i}\qquad \forall~i>0
  $$ 
 and
  $$ 
 RT^\zeta (  \tau(\varrho)*\alpha)=t_0(\varrho)RT^\zeta(\alpha) \neq 0
  $$ 
  \qed

  \begin{theorem} Let $\zeta$ be a   root of unity whose order is not divisible by $4$. If $X$ is an oriented three-dimensional homology sphere with finite character variety then
    $$
     \dim_{\mathbb{C}} \Kz(X)\geq \dim_{\mathbb{C}}\Ke(X). 
     $$ \end{theorem}
  \proof Theorem \ref{itsbig} handles the case when the order of $\zeta$  is even. If the order is  odd then the conclusion follows from Proposition \ref{allgood}.
  \qed

\begin{rem}
  
  It is natural to ask whether for a homology sphere it can happen that 
  $
  \dim_{\mathbb{C}}\Kz(X)_{\fm_0}>\dim_{\mathbb{C}}\Ke(X)_{\fm_0}
  $? By \cite{W}, if $X$ is a homology sphere, then $\fm_0$ is a regular point of the character variety of $\pi_1(X)$; thus, we are asking if the dimension of $\Kz(Z)_{\fm_0}$ can be greater than $1$ in this case. 
  \end{rem}
  \vskip.1in

\subsection{The dimension of the $\mathbb{Q}(q)$-skein module}
 Let $\Phi_{2m}(q)$ denote the $2m$-th cyclotomic polynomial.
  Following \cite{DKS}, we say that a $\mathbb{Q}[q,q^{-1}]$ module $M$  is {\bf tame} if 
  \begin{itemize} [wide, nosep, labelindent = 0pt, topsep = 1ex]
  \item $M$ is a direct sum of cyclic modules,
  \item and, for some odd $m$, it does not contain a submodule that is isomorphic to $\mathbb{Q}[q,q^{-1}]/(\Phi_{2m}(q))$.
  \end{itemize}   
  This means that $M$ can be written as
  $$ 
  M=F\oplus T
  $$ 
  where $F$ is free and $T$ is torsion.  By assumption, there exists odd $m$  such that no submodule has
  $\Phi_{2m}(q)$-torsion. From \cite{DKS} if the universal skein module of $X$ is tame then it has a finite character variety.

  Suppose $\zeta$ is a primitive $2m$-th root of unity. Make $\mathbb{C}$ into a $\mathbb{Q}[q,q^{-1}]$-module by letting $q$ act as multiplication by $\zeta$. Notice that
  $$ 
  M\otimes_{\mathbb{Q}[q,q^{-1}]} \mathbb{C}= F\otimes_{\mathbb{Q}[q,q^{-1}]} \mathbb{C}
  $$
  as  all the torsion elements are killed by the tensor product.\\

  The following is a corollary of Theorem~\ref{itsbig}.

  \begin{cor} \label{otherbound} If $X$ is a homology sphere with  tame universal skein module  then
    $$ 
    \dim_{\mathbb{Q}(q)}K_q(X) \geq  \dim_\mathbb{C} \Ke(X).
    $$ 
  \end{cor}

  \proof By \cite{DKS} if the universal skein module is tame then the character variety of $X$ is finite. Starting with the universal skein module $K(X)$, by assumption  $$K(X)\otimes_{\mathbb{Z}[q,q^{-1}]}\mathbb{Q}[q,q^{-1}]=F \oplus T$$where
    $F$ is free and $T$ is torsion. Also there is an odd $m$ so that $T$ has no $\Phi_{2m}$-torsion.
    To compute $K_q(X)$, we tensor $F\oplus T$ with $\mathbb{Q}(q)$. Since  $\mathbb{Q}(q)$ is torsion free, we  find
    $$ 
    K_q(X)=F\otimes_{\mathbb{Q}[q,q^{-1}]}\mathbb{Q}(q). 
    $$ 
    On the other hand, computing $\Kz(X)$ with the conditions that 
    \begin{itemize} [wide, nosep, labelindent = 0pt, topsep = 1ex]
    \item $\zeta$ has order $2m$ and $m$ is odd,
    \item $K(X)$ has no $\Phi_{2m}$-torsion,
    \end{itemize}
    we get
    $$ 
    \Kz(X)=F\otimes_{\mathbb{Q}[q,q^{-1}]} \mathbb{C},
    $$ 
    where $q$ acts as multiplication by $\zeta$ on $\mathbb{C}$.
    This means the two modules have the same dimension over the corresponding fields. By Theorem \ref{itsbig}, we conclude
    $$ \dim_{\mathbb{Q}(q)}K_q(X)=\dim_{\mathbb{C}}\Kz(X) \geq \dim_{\mathbb{C}}\Ke(X). 
    $$ \qed

      \begin{theorem}\label{rational} If $X$ is a homology sphere with tame universal Kauffman bracket skein module then 
      $$
      \dim_{\mathbb{C}}\Ke(X)=\dim_{\mathbb{Q}(q)}K_q(X).
      $$ \end{theorem}

    \proof From \cite{DKS} $\dim_{\mathbb{Q}(q)}K_q(X)\leq \dim_{\mathbb{C}}\Ke(X)$.  Together with  Corollary \ref{otherbound} this  yields the result.
\qed\\

 Suppose that $\fm$ is a point of the character variety of the manifold $X$ corresponding to an irreducible  representation $\rho$.
The {\bf normal cone}
of $\rho$ is a variety whose unreduced coordinate ring is  
$$ 
\textnormal{Cone}(\rho)=\bigoplus_{i=1}^{k-1} \fm^i/\fm^{i+1}.
$$
The normal cone describes infinitesimal deformations of the representation $\rho$.
Theorem~\ref{main_thm} implies that there is an injective complex linear map,
$$ 
\theta_\rho\colon \textnormal{Cone}(\rho)\rightarrow \Kz(X).
$$ 
Therefore, given root of unity $\zeta$, an  irreducible representation $\rho$, and an infinitesimal deformation of $\rho$, there is an invariant of framed links in $X$ that satisfies the Kauffman bracket skein relation. 
 \begin{rem} 
   It would be interesting to understand the invariants corresponding to infinitesimal deformations of $\rho$  on the level of tangle functors. 
 \end{rem}

 \begin{rem} To understand the local structure of the Kauffman bracket skein module at a representation that is reducible, it is probably better to work at the level of the stated skein algebra \cite{Le2,CL,LY}, following the approach of \cite{GJS1} where the skein algebras correspond to nonsingular objects.
\end{rem}

\bibliographystyle{hamsalpha}
\bibliography{biblio}

\end{document}